\documentclass{article}
\usepackage{amsmath,amssymb,amsfonts,theorem}%
\usepackage{graphicx,color}
\graphicspath{{./Figures/}}

\newcommand{\until}[1]{\{1,\dots, #1\}}
\newcommand{\subscr}[2]{#1_{\textup{#2}}}

\newcommand{\setdef}[2]{\{#1 \; | \; #2\}}

\newcommand{\union}{\operatorname{\cup}}

\newcommand{\naturals}{\mathbb{N}} 
\newcommand{\reals}{\mathbb{R}} 

\newcommand{\integersnonnegative}{\ensuremath{\mathbb{Z}}_{\ge 0}}
\newcommand{\integernonnegative}{\ensuremath{\mathbb{Z}}_{\ge 0}}

\newcommand{\realpositive}{\ensuremath{\reals_{>0}}}

\newcommand{\N}{\mathbb{N}}
\newcommand{\Nzero}{\integersnonnegative}
\newcommand{\R}{\reals}

\newcommand{\1}{\mathbf{1}} 
\newcommand{\G}{\mathcal{G}}    

\renewcommand{\S}{{\mathcal{S}}}

\newcommand{\card}[1]{\left|#1\right|}  

\newcommand{\xave}{\subscr{x}{ave}}

\newcommand{\kin}{\subscr{k}{in}} 
\newcommand{\kout}{\subscr{k}{out}}

{ \theorembodyfont{\normalfont} 

\newtheorem{remark}{Remark}
}

\newtheorem{definition}{Definition}
\newtheorem{theorem}{Theorem}
\newtheorem{lemma}[theorem]{Lemma}

\title{Efficient quantization for average consensus}
\author{Ruggero Carli, Fabio Fagnani, Paolo Frasca, Sandro Zampieri}
\date{\today}
\begin{document}
\maketitle
\begin{abstract}
This paper presents an algorithm which solves exponentially fast the average consensus problem on strongly connected network of digital links. The algorithm is based on an efficient {\em zooming-in/zooming-out} quantization scheme.
\end{abstract}

\section{Introduction}
Problems of coordination, and in particular of distributed averaging, have become an important subject of research for control theorists \cite{ROS-JAF-RMM:07}. In the average consensus problem, $N$ agents are given one number each, and have to compute their average.  The interest of the problem comes from the communication limitations the agents are subject to, which require the design of an iterative algorithm. Usually, the agents are allowed to communicate with a restricted number of neighbors, and this is described by a suitable directed graph. In this work, we present an average consensus algorithm in which the agents can communicate with their neighbors only via rate constrained digital channels. We propose an efficient encoding/decoding strategy, based on the so called {\it zooming-in/zooming-out} quantization method \cite{RWB-DL:00,ST-SM:04,GNN-FF-SZ-RJE:07}, and we give theoretical and simulation results on its performance. The algorithm solves the average consensus problem with arbitrary precision in exponential time, using a quantizer with number of states which is finite, but can depend on $N$.

\subsection{Related works}
The constraint of quantization of messages, or of states, has been considered in several recent papers \cite{LX-SB-SL:05,RE-BM-SS:06,RC-FF-PF-TT-SZ:07,PF-RC-FF-SZ:08,AK-TB-RS:07,SK-JM:sub,AN-AO-AO-JNT:07,TCA-MJC-MGR:08}.
However, existing algorithms in the literature are not able to achieve average consensus with arbitrary precision, with the notable exception of \cite{RC-FB-SZ:08j}, which uses logarithmic quantizers, together with a coding scheme similar to the one in the present paper. The zooming-in/zooming-out algorithm and the main convergence result of this paper appeared, in slightly different version, in \cite{RC-FF-PF-SZ:07} and in the thesis \cite{PF:09}. Note that recently \cite{GB-SZ:09} a modification of the present algorithm has been proved to converge under a different condition.

\subsection{Outline}
The problem of average consensus is formally presented in Section~\ref{sec:problem}. The algorithm, and the main result, are in  Section~\ref{sec:algo}. In Section~\ref{sec:simul} we present several simulations, which show an encouraging performance, beyond the scope of the proven result. Finally, Section~\ref{sec:outro} points out some potential developments.

\section{Average consensus}\label{sec:problem}
Let there be a directed graph $\G=(V,E)$, with $V=\until{N}$ a set of dynamical systems, that we call agents. We assume that only if $(j,i) \in \,E$, $j$ can transmit information about its state to $i$. Let us consider, for $t\in \integernonnegative$, the following discrete-time state equations
\begin{align*}
     x_i(t+1)=x_i(t)+u_i(t)\qquad i=1,\ldots,N
\end{align*}
where $x_i(t)\in\reals$ is the state of the $i$-th system/agent, and $u_i(t)\in\reals$ is the control input.
More compactly we can write
\begin{equation}\label{eq:stateeq}
     x(t+1)=x(t)+u(t)
\end{equation}
where $x(t),u(t)\in\reals^N$.

\begin{definition}[Average consensus problem]\label{def:ConsensusProblem}
Given the system \eqref{eq:stateeq}, the {\em average consensus problem} consists in designing a sequence of control inputs $u(t)$ yielding the consensus of the states, namely
\begin{equation}\label{eq:agree}
\lim_{t\to+\infty}x(t)=\xave(0) \1,
\end{equation}
where $\1$ is a $N$-vector whose entries are 1.
\end{definition}
The important case we restrict to, and which has been mostly studied in literature, is a static linear feedback
\begin{equation}\label{eq:feedback}
u(t)=Kx(t)\qquad K\,\in\,\reals^{N \times N}
\end{equation}
In such case the system (\ref{eq:stateeq}) is given by the
following closed loop system
\begin{equation}\label{eq:closedloopAlgo}
x(t+1)=(I+K)x(t).
\end{equation}
The matrix $I+K$ is commonly called a Perron matrix, and denoted by $P$.
The problem of designing the controller thus reduces to find a matrix $K$ {\em adapted to $\G$}, in the sense that if $(j,i) \notin \,E$, then $K_{ij}=0$. Methods for design are known in the literature \cite{RC-FF-AS-SZ:08}. We shall only recall here that the average consensus problem is solved if and only if the matrix $P$ is a doubly stochastic matrix.

\subsection{Quantization}
We can make the framework \eqref{eq:feedback} more general if we assume that the agents can not directly access their neighbors' states, but can obtain estimates or approximations of them. Let $\hat x_j^i(t)$ be the estimate agent $i$ has of the state of agent $j$, at time $t$.
The control input $u_i(t)$ assumes the following form
\begin{equation}\label{eq:generalChannels}
u_i(t)= \sum_{j=1}^N K_{ij}\hat{x}_{j}^i(t).
\end{equation}
Due to the possibility that $\hat{x}_{j}^i(t)\neq x_j(t)$, such a control is not in general guaranteed to yield convergence to consensus, nor to preserve the average of states.

Most of the literature on the consensus problem assumes that the communication channels between the nodes allows to communicate real numbers with no errors. In practical applications this can be an unrealistic assumption: the communication is in most cases digital and the channel is noisy. There can be strict bandwidth constraints, as well as communication delays, and packet losses due to interferences and erasures.
These limitations are mostly likely to be significant for a network communicating in a wireless fashion, which is nowadays a typical choice. In this paper, we concentrate on modeling the links as instantaneous and lossless digital channels. This choice forces a quantization on the real numbers that agents have to transmit, but keeps aside the issues of delays and packet losses.

\section{Encoding-decoding schemes}
From now on, let us assume to have a strongly connected directed graph $\G$ and an average consensus controller $K$ adapted to $\G$.
In the extended communication framework \eqref{eq:generalChannels}, the control input $u_i(t)$ of agent $i$ has the following form
\begin{equation}\label{eq:input}
u_i(t)= \sum_{j=1}^N K_{ij}\hat{x}_{j}^i(t)\,,
\end{equation}
where $\hat{x}_{j}^{i}$ is the estimate of the state $x_j$ which has
been built by the agent $i$.
The point now is to describe how the estimates are built. In other works, \cite{PF-RC-FF-SZ:08}, \cite{RC-FF-PF-SZ:08}, we have assumed that the estimate be just the transmitted message. Now, instead, we look for a non trivial way of constructing estimates.

Suppose that the $j$-th agent sends to the $i$-th agent, through a
digital channel, at each time instant $t$, a symbol $s_{ij}(t)$
belonging to a finite or countable alphabet $\mathcal{S}_{ij}$,
called the \emph{transmission alphabet}. It is assumed that the
channel is reliable, that is each symbol transmitted is received
without error. In general, the structure of the coder by which the
$j$-th agent produces the symbol to be sent to the $i$-th agent can
be described by the following equations
\begin{equation}\label{coderstate}
\left\{
\begin{array}{rcl}
\xi_{ij}(t+1)&=&F_{ij}(\xi_{ij}(t),s_{ij}(t))\\
s_{ij}(t)&=&Q_{ij}(\xi_{ij}(t),x_j(t))\\
\end{array}\right.
\end{equation}
where $s_{ij}(t)\,\in\,\mathcal{S}_{ij}$,
$\xi_{ij}(t)\,\in\,\Xi_{ij}$, $Q_{ij}:\Xi_{ij}\times \reals
\to\mathcal{S}_{ij}$, and $F_{ij}:\Xi\times
\mathcal{S}_{ij}\to\Xi_{ij}$.
The decoder, run by agent $i$, is the system
\begin{equation}\label{decoderstate}
\left\{
\begin{array}{rcl}\xi_{ij}(t+1)&=&F_{ij}(\xi_{ij}(t),s_{ij}(t))\\
\hat{x}_{j}^i(t)&=&H_{ij}(\xi_{ij}(t),s_{ij}(t)),
\end{array}\right.
\end{equation}
where $H_{ij}:\Xi_{ij}\times \mathcal{S}_{ij}\to \reals$.
The set $\Xi_{ij}$ serves as \emph{state space} for the coder/decoder, whereas the maps $F_{ij}, Q_{ij}, H_{ij}$ represent, respectively, the \emph{coder/decoder dynamics}, the \emph{quantizer function}, and the \emph{decoder function}. Coder and decoder are jointly initialized at $\xi_{ij}(0)=\xi_0$.\\
In general, one may have different encoders at agent $j$, according to the various agents the agent $j$ wants to send its data. For the sake of simplicity, we assume in the sequel that each system uses the same encoder for all data transmissions. Thus, each agent $j$ broadcasts the same symbol $s_j(t)$ to all its neighbors. In this case every agent receiving data from $j$ obtains the same estimate of $x_j(t)$, which we shall denote by $\hat{x}_{j}(t)$. In this way, by letting $F_j=F_{ij}, H_j=H_{ij}, Q_j=Q_{ij}$ and $\Xi_{j}=\Xi_{ij}$, the previous coder/decoder couple can be represented by the following {\em state estimator with memory}
\begin{equation}\label{stateestimator}
\left\{
\begin{array}{rcl}
\xi_j(t+1)&=&F_j(\xi_j(t),s_j(t))\\
s_j(t)&=&Q_j(\xi_j(t),x_j(t))\\
\hat{x}_j(t)&=&H_j(\xi_j(t),s_j(t))
\end{array}
\right.
\end{equation}
Moreover (\ref{eq:input}) assumes the following form
\begin{equation}\label{eq:input1}
u_i(t)= \sum_{j=1}^N K_{ij}\hat{x}_{j}(t),
\end{equation}
and the system evolution is
\begin{equation}\label{eq:EvolWithEstim}
x(t+1)=x(t)+K \hat{x}(t).
\end{equation}

\begin{remark}[Need for synchrony]\label{rem:NeedForSynch}
Remark that, from the point of view of the implementation, both agent $j$ and its neighbors have to run identical copies of the system $\xi_j(t).$ This clearly implies an intrinsical weakness of the scheme with respect to failures: if messages are lost, the synchrony of these copies is lost as well, and the algorithm does not work.
\end{remark}
In change of this potential drawback, the scheme has very interesting convergence properties.

\section{Zooming-in zooming-out strategy}\label{sec:algo}
Our strategy is inspired by the quantized stabilization
technique proposed in \cite{RWB-DL:00}, which is called zooming-in/zooming-out strategy. In this case the information exchanged between
the agents is quantized by scalar uniform quantizers which assume values in a finite set, and can be
described as follows. For $m\,\in\,\naturals$ define the set of
quantization levels
$$
\mathcal{S}_m=\setdef{-1+\frac{2\ell-1}{m}}{\ell\,\in\,\until{m}}\union \{-1,1\}.
$$
The corresponding \emph{uniform quantizer} $q^{(m)}:\reals\rightarrow \mathcal{S}_m$ is as follows. Let $x\,\in\,\reals$ then
$$
q^{(m)}(x)=-1+\frac{2\ell-1}{m}
$$
if $\ell\,\in\,\left\{1,\ldots,m\right\}$ satisfies
$-1+\frac{2(\ell-1)}{m}\leq x \leq -1+\frac{2\ell}{m}$, otherwise
$q^{(m)}(x)=1$ if $x>1$ and $q^{(m)}(x)=-1$ if $x<-1$.

It is easy to see that the quantizer $q^{(m)}$ enjoys the following property.
\begin{lemma}[Quantization error]\label{lem:quantErrBound}
Given $z\in \R$ and $l\in \realpositive$ such that $|z|\le l$, it holds
\begin{align}
|z- l q^{(m)}\left(\frac{z}{l}\right)|\le \frac{l}{m}.
\end{align}
\end{lemma}

\smallskip

Let $m\,\in \, \naturals$, $\kin \,\in\,]0,1[$, and $\kout\,\in\,]1,+\infty[.$
The {\em zooming-in/zooming-out} encoder/decoder, with parameters $m,\kin,\kout$, is defined by the
alphabet $\mathcal{S}=\mathcal{S}_{m}$, the state space
$\Xi_j=\reals\times \reals_{>0}$, such that
$\xi_j(t)=\left(\hat x_j(t),l_j(t)\right)$, and the
dynamics
\begin{align*}
&\hat x_j(0)=0\\
&\hat x_j(t+1)=\hat x_j(t)+l_j(t+1) s_j(t+1) \qquad \forall \,t\ge 0,
\end{align*}
and
\begin{align*}
& l_j(0)=l_0 \in \R\\
& l_j(1)=l_0 \\
& l_j(t+1)=
\begin{cases}
\kin l_j(t)\qquad & \mbox{if}\,\, |s_j(t)|<1 \\
\kout l_j(t)\qquad & \mbox{if}\,\, |s_j(t)|=1
\end{cases}\qquad \forall \,t>0.
\end{align*}

The sent symbol is
\begin{eqnarray*}
s_j(t)&=&q^{(m)}\left(\frac{x_j(t)-\hat x_j(t-1)}{l_j(t)}\right)\qquad \forall \,t>0.
\end{eqnarray*}

Let us comment on such definition. Remark that the first component of the
coder/decoder state contains $\hat x(t)$, the estimate of $x(t)$. The transmitted
messages contain a quantized version of the estimation error $x_j(t)-\hat x_j(t-1)$
scaled by the factor $l_j(t)$. Accordingly, the second component of
the coder/decoder state, $l_j$, is referred to as the
\emph{scaling factor}: it increases when $|x_j(t)-\hat{x}_j(t-1)|>l_j(t)$
(``zooming out step") and decreases when ${|x_j(t)-\hat{x}_j(t-1)|\leq l_j(t)}$ (``zooming in step").

We can now prove the main theorem of this chapter.
\begin{theorem}[Convergence]\label{th:ZoomConv}
Consider the system (\ref{eq:EvolWithEstim}) with the zooming coding.
Let $\rho$ be the essential spectral radius of $I+K$. Suppose that $\rho< \kin <1$, $m\geq
\frac{(4+3\kin )\sqrt{N}}{\kin (\kin -\rho)}$ and that
$\l_0>\frac{2(\rho+2)\|x(0)\|}{\kin -\frac{3\sqrt{N}}{m}}$.
Then, for any initial condition $x(0)\,\in\,\reals^N$,
$$
\lim_{t\rightarrow
+\infty}x(t)=\lim_{t\rightarrow+\infty}\hat{x}(t)=\xave\1,
$$
and the convergence is exponential, with rate not bigger than $\kin.$
\end{theorem}
\begin{proof}%
First, observe that $\1^* K=0$, and then the states update rule preserves the average of the state $x(t)$ at each time step. Hence if we prove convergence to a consensus value, that will be the average.
Let us rewrite the overall system as the coupling of the three dynamical systems on $\R^N$,
\begin{align*}
&l_j(0)=l_0 \qquad\,\forall j\in V\\
&l_j(1)=l_0 \qquad\,\forall j\in V\\
&l_j(t+1)=
\begin{cases}
\kin l_j(t)\qquad & \mbox{if}\,\, |x_j(t+1)-\hat x_j(t)|<1 \\
\kout l_j(t)\qquad & \mbox{if}\,\, |x_j(t+1)-\hat x_j(t)|\ge 1
\end{cases}\qquad\,\forall j\in V, t\in \N
\end{align*}
\begin{align*}
&x(0)=x_0\\
&x(t+1)=x(t)+K \hat x(t) \qquad\,\forall t\in \Nzero
\end{align*}
\begin{align*}
&\hat x_j(0)=0\qquad\,\forall j\in V\\
&\hat x_j(t+1)=\hat x_j(t) +l_j(t+1) q^{(m)}\left(\frac{x_j(t+1)-\hat x_j(t)}{l_j(t+1)}\right) \qquad\,\forall j\in V, t\in \Nzero .
\end{align*}
Their coupling implies that the following facts are equivalent:
\begin{enumerate}
\item $\displaystyle\lim_{t\to+\infty} l(t)=0;$
\item $\displaystyle\lim_{t\to+\infty} x(t)=\xave\1;$
\item $\displaystyle\lim_{t\to+\infty} \left(\hat x(t)-x(t)\right)=0.$
\end{enumerate}
The idea of the proof is the following: we show that under the assumptions there always happen zooming in steps. This implies $l(t)=\kin^{t-1}l_0\1$ and implies in turn convergence to average consensus.

By definition, there are only zooming in steps if
\begin{align}
|x_j(t+1)-\hat x_j(t)|< l_j(t+1) \quad \forall i\in V, t\in \Nzero
\end{align}
or equivalently if
\begin{align}\label{eq:condZoom}
|x_j(t+1)-\hat x_j(t)|< \kin^t l_0 \quad \forall i\in V, t\in \Nzero .
\end{align}

Let us define $\tilde y(t)=K x(t)$ and $e(t)=x(t)-\hat{x}(t),$
and remark that $\eqref{eq:condZoom}$ is equivalent to
\begin{equation}
|\tilde y(t)-(K-I) e(t)|<\kin^t l_0,
\end{equation}
where the inequality is meant componentwise.

To prove this fact, we proceed by strong induction over $t$, proving that for all $t\in \Nzero$
\begin{align}
& \|\tilde y(t)\|\le \kin^t l_0 \label{eq:condZoom2a}\\
& \| \tilde y(t)-(K-I) e(t)\|\le \kin^t l_0 \label{eq:condZoom2b}.
\end{align}
Since $\| \tilde y(t)-(K-I) e(t)\|_\infty\le\| \tilde y(t)-(K-I) e(t)\|,$ the inequality \eqref{eq:condZoom2b} implies \eqref{eq:condZoom} and then convergence.

Let us thus check the validity of Equations~\eqref{eq:condZoom2a}~and~\eqref{eq:condZoom2b} for $t=0$, $t=1$, and for a general $t+1$, with $t\in \N$, given it true for all $t'<t$.
We shall use the recursion
\begin{equation*}
\tilde y(t+1)=(I+K)\tilde y(t)-K^2e(t),
\end{equation*}
as well as Lemma~\ref{lem:quantErrBound}, which implies $\|e(t)\|\le \frac{\sqrt{N}}{m}l_0.$

Using triangle and submultiplicative inequality of norms, and ${\|I+K\|=\rho},$ ${\|K\|\le 2}$, and ${\|K-I\|\le 3}$, we get that
\begin{description}
\item[$t=0.$] \eqref{eq:condZoom2a} holds if
$\|K x_0\|<l_0,$ and
\eqref{eq:condZoom2b} holds if $\|x_0\|_\infty<l_0.$
\item[$t=1.$] \eqref{eq:condZoom2b} holds if
\begin{equation}\label{zoomCondA}
 2 (2+\rho)\|x_0\|+3 \frac{\sqrt{N}}{m}l_0<\kin l_0,
 \end{equation}
  and \eqref{eq:condZoom2a} holds if
  \begin{equation}\label{zoomCondB}{\|K x_0\|<\kin l_0}.\end{equation}
\item[$t>1.$] Assuming \eqref{eq:condZoom2a} and \eqref{eq:condZoom2b} true for time steps $t$ and $t-1$, we obtain that \eqref{eq:condZoom2b} holds for the time step $t+1$ if
$$\rho\kin^t+4 \frac{\sqrt{N}}{m}\kin^{t-1}l_0 + 3 \frac{\sqrt{N}}{m}\kin^t l_0<\kin^{t+1}l_0,$$
that is if
\begin{equation}\label{zoomCondC}\rho\kin+4 \frac{\sqrt{N}}{m} + 3 \frac{\sqrt{N}}{m}\kin <\kin^{2}.\end{equation}
Instead, \eqref{eq:condZoom2a} holds if
\begin{equation}\label{zoomCondD}\rho\kin+4 \frac{\sqrt{N}}{m} <\kin^{2}.\end{equation}
\end{description}
The system of conditions \eqref{zoomCondA}-\eqref{zoomCondB}-\eqref{zoomCondC}-\eqref{zoomCondD} is satisfied when the assumptions of the theorem hold.
\end{proof}
This result is quite conservative in proving convergence, since it guarantees that since the beginning only zooming in steps happen. This is clearly restrictive, since there is no need for $l(t)$ to be monotonic to converge to zero.

\section{Simulation results}\label{sec:simul}
In this section we show by simulations the properties of the algorithm in terms of (speed of) convergence, as it depends on the parameters of the method, $\kin,\kout$ and $m$, and on the topology of the graph. Indeed, the convergence result we obtained, Theorem~\ref{th:ZoomConv}, gives sufficient conditions on the parameters, depending on $\rho$, and on $N$, and thus on the graph. Convergence is possible, provided the conditions, at exponential speed, with rate $\kin>\rho.$ Simulations are worth of interest, since they demonstrate that convergence is possible also {\em outside the scope of Theorem~\ref{th:ZoomConv}}, and with a {\em speed which can be faster than the linear algorithm with ideal communication \eqref{eq:closedloopAlgo}.}

Roughly speaking, we can identify two different regimes for the evolution of the algorithm, depending on the zooming-in rate $\kin$ which is enforced. If $\kin>\rho$, the zooming-in follows the natural contraction rate of the system, and according to the spirit of the proof of Theorem~\ref{th:ZoomConv}, eventually (almost) only zooming-in steps happen, and the rate of convergence is (no better than) $\kin$. The choice of $\kout$ plays essentially no role.
If instead $\kin<\rho$, the state does not contract fast enough for a stationary sequence of zooming-ins to establish. Then, zooming-in and zooming-out steps alternate in a complicated way, leading anyhow to consensus, with an evolution  which depends on both $\kin$ and $\kout$: the convergence can be faster than in the ideal case.
Such phenomenon is rather surprising at a first glance, but we can explain it intuitively if we consider that the zooming-in/zooming-out algorithm takes advantage of the use of memory, and thus the sent messages are actually more informative.

Most of the simulations regard ring graphs, because concentrating on a simple example of a poorly connected graph helps to highlight the features of the algorithm. Weights of matrices $P$ have been chosen according to a maximum degree rule. Hence a ring of $N=20$ nodes induces a matrix whose essential spectral radius is $\rho=0.9673$. The initial conditions have been generated according to a gaussian distribution, and all simulations with $20$ nodes start from the same initial condition.

\begin{figure}[htb]\centering
\includegraphics[width=.8\textwidth]{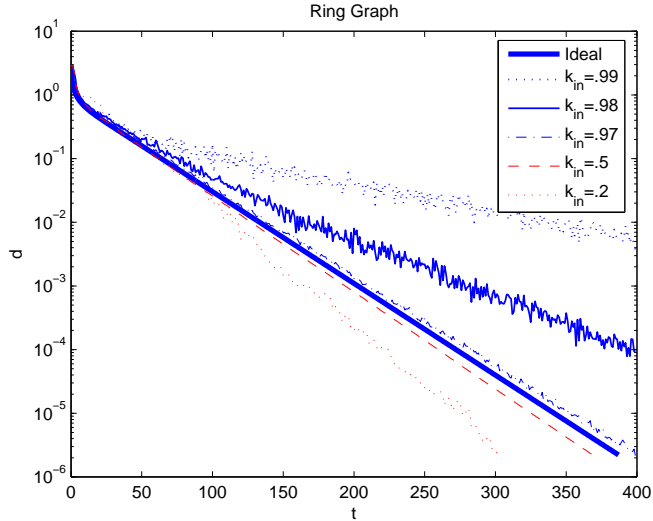}
\caption{Ring graph. Dependence on $\kin$: remark the threshold $\rho=0.9673$. $N=20$, $m=6$, $\kout=2$.}\label{fig:ZoomRingkin}
\end{figure}

\begin{figure}[htb]\centering
\includegraphics[width=.8\textwidth]{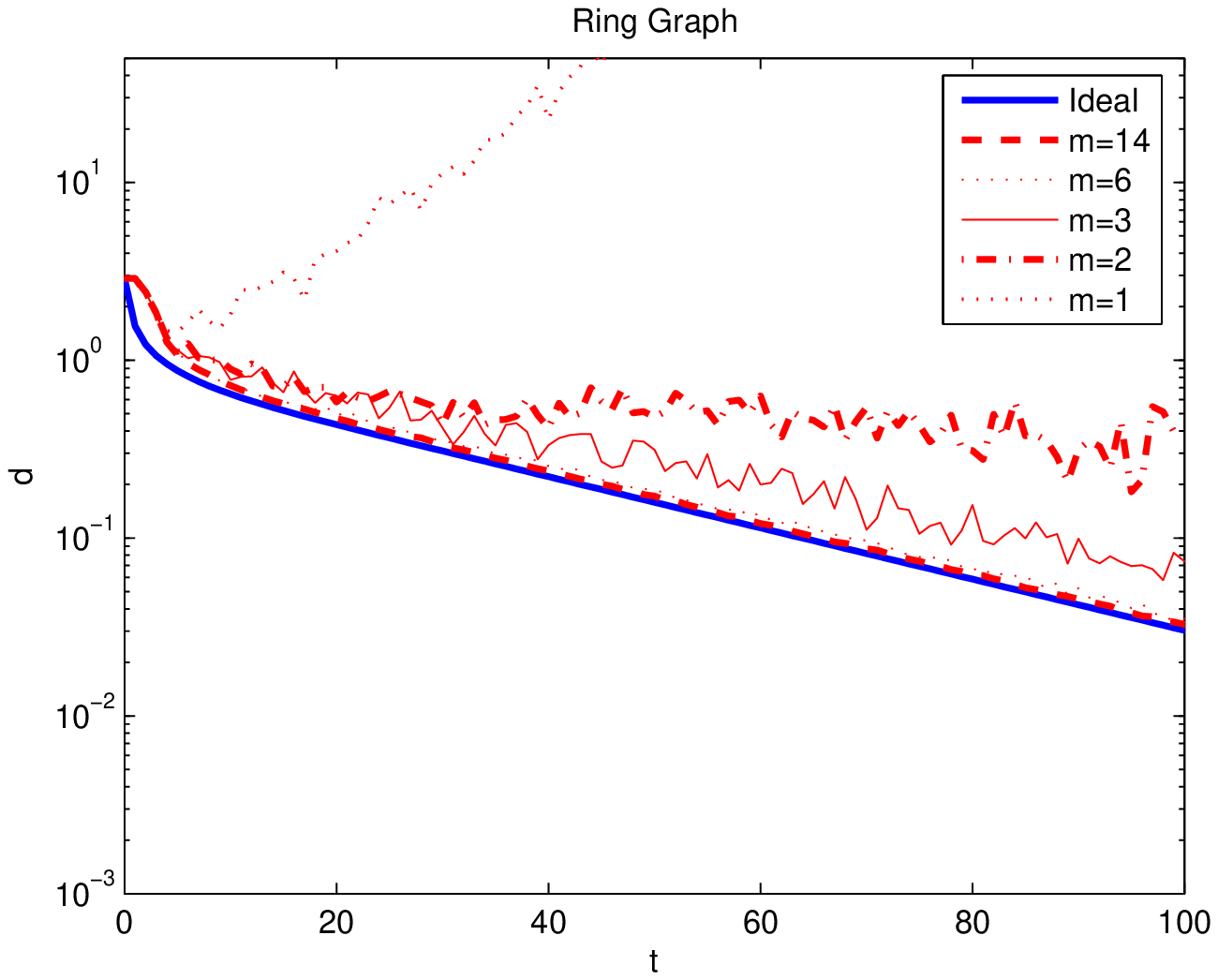}
\caption{Ring graph, in the $\kin>\rho$ regime. $N=20$, $\kin=0.97$, $\kout=2$.}\label{fig:ZoomRing_m1}
\end{figure}

The transition between the two regimes is explored in Figure~\ref{fig:ZoomRingkin}, confirming the role of $\rho$. Instead, Figures~\ref{fig:ZoomRing_m1}~and~\ref{fig:ZoomRing_m1} consider the role of the number of quantization levels $m$. In the $\kin>\rho$ regime, at least $m=3$ is required for convergence, and increasing $m$ leads to improve the rate, approaching $\kin$. In the $\kin<\rho$ regime, $m=1$ is enough for convergence, and increasing $m$ leads to slightly improve the rate. Since $\card{\S_m}=m+2$, then $5$ and $3$ symbols are sufficient for convergence, respectively. If, taking advantage of the synchrony, one encodes the middle level, which is zero whenever $m$ is odd, as a no-signal, the required number of level can be obtained with only two bits (one bit, respectively).

\begin{figure}[htb]\centering
\includegraphics[width=.8\textwidth]{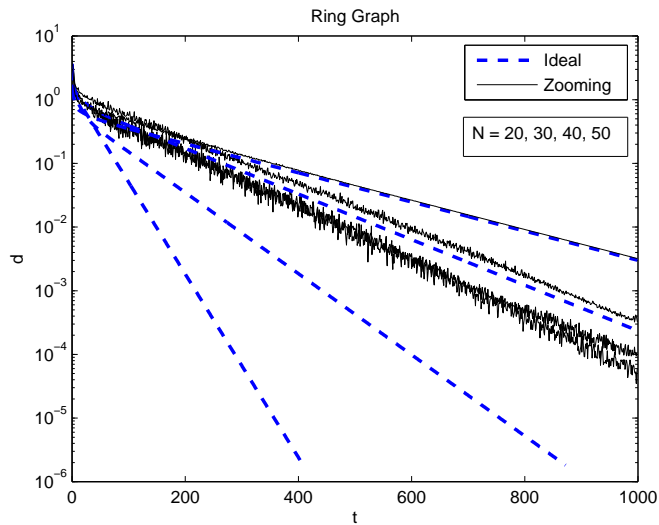}
\caption{Ring graph, with $m=3$, $\kin=0.99$, $\kout=2$. For $N>20$ we are out of the scope of Theorem~\ref{th:ZoomConv}, and for $N>30$, we have $\kin<\rho$. }\label{fig:ZoomRingFixedParam2}
\end{figure}

\begin{figure}[htb]\centering
\includegraphics[width=.49\textwidth]{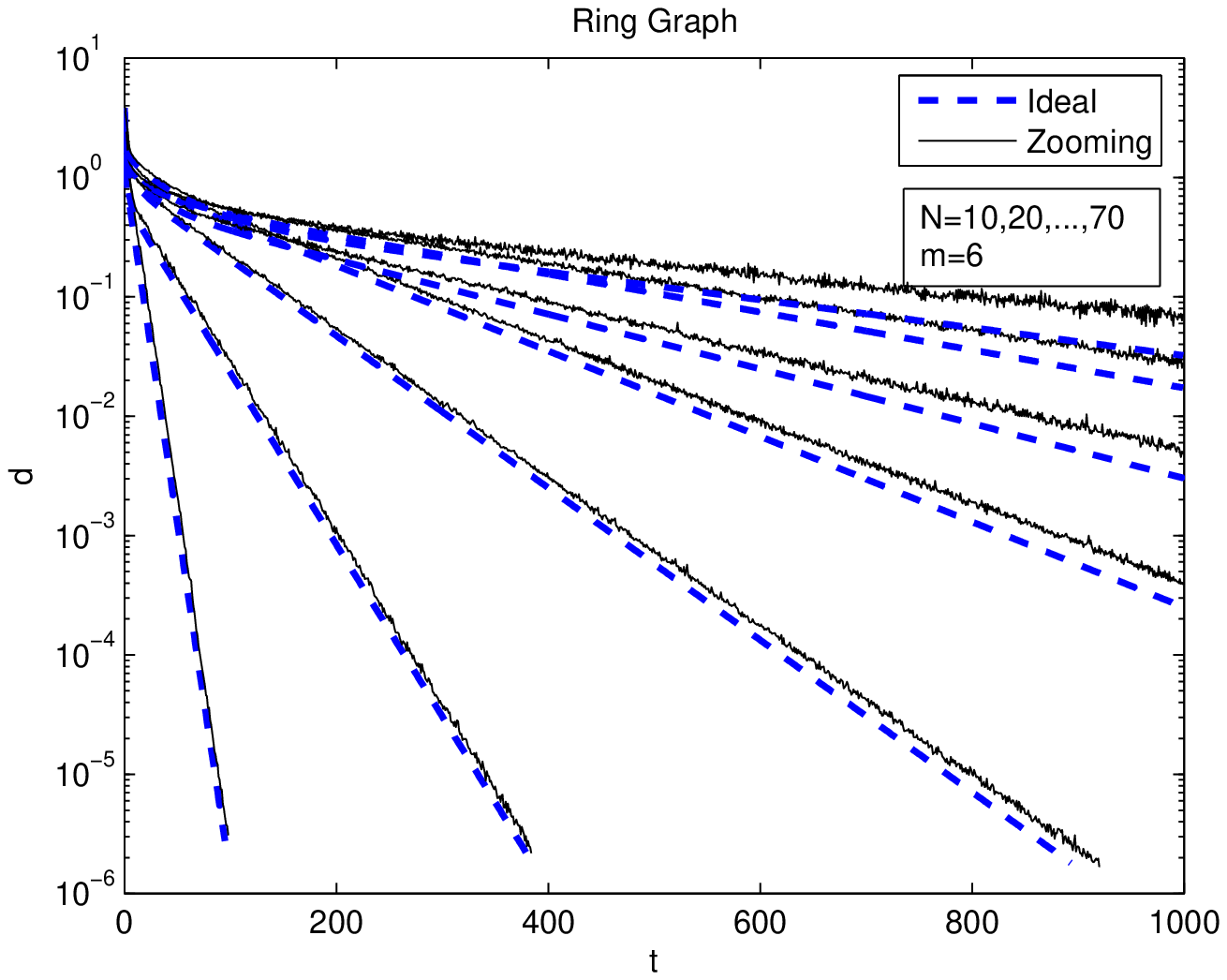}
\includegraphics[width=.49\textwidth]{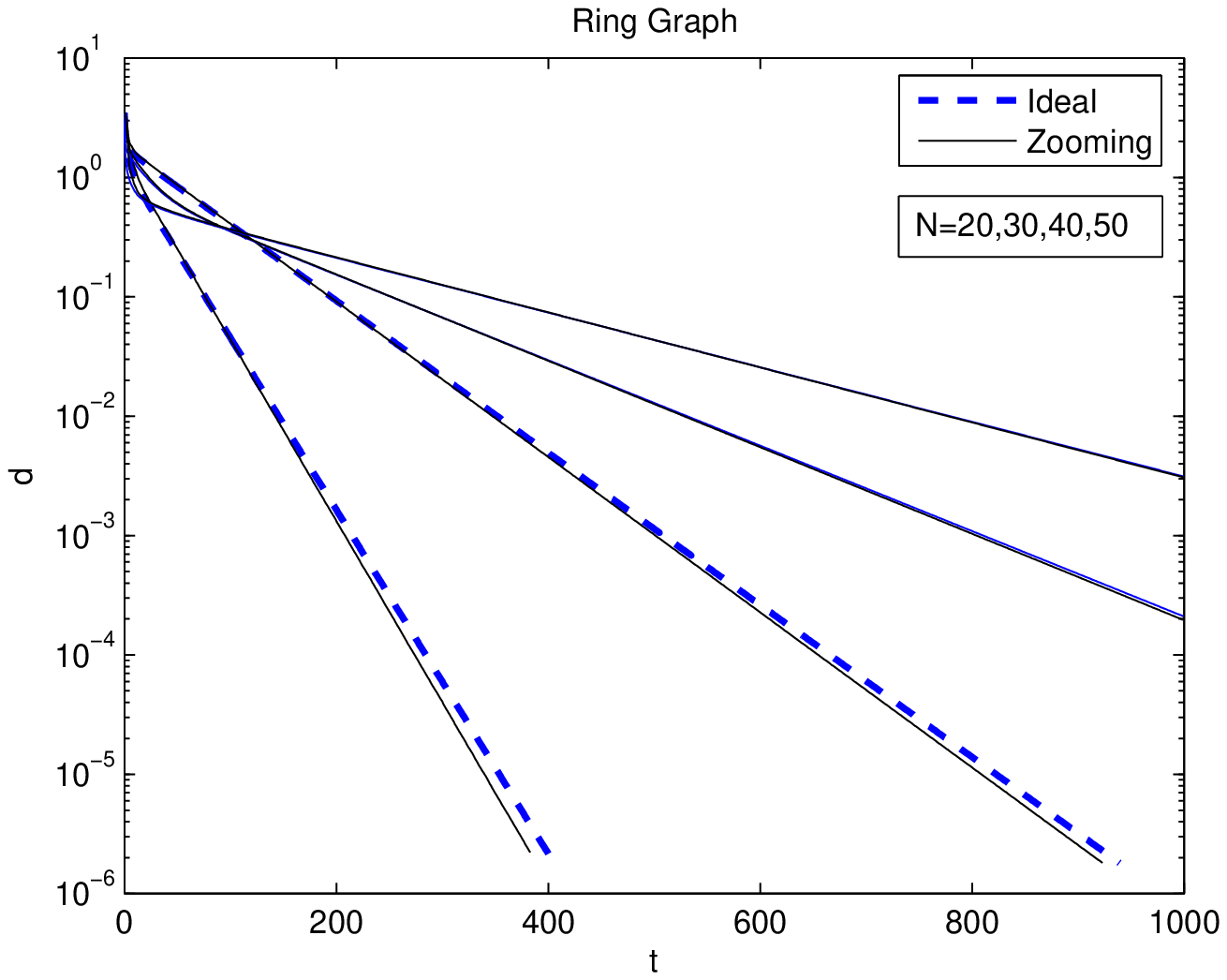}
\caption{Ring graph, $m=6$, for increasing $N$. Left plot assumes $\kin=1.01\rho$ and $\kout=2$, yielding the $\kin>\rho$ regime. Right plot assumes $\kin=0.5\rho$, $\kout=\kin^{-1}$, yielding the $\kin<\rho$ regime. 
}\label{fig:ZoomRing_scales}
\end{figure}

The scaling properties of the algorithm are a main point of interest.
Figure~\ref{fig:ZoomRingFixedParam2} considers how the performance of the algorithm (with fixed parameters) degrades as $N$ increases, in a sequence of ring graphs. The plot shows that $m$ need not to depend on $N$ as required Theorem~\ref{th:ZoomConv}. Instead, Figure~\ref{fig:ZoomRing_scales} considers two cases in which the zooming factors $\kin$ and $\kout$ are depends on the graph, in such a way to keep the evolution in the same regime, in spite of the growth of $N$.

\begin{figure}[htb]\centering
\includegraphics[width=.8\textwidth]{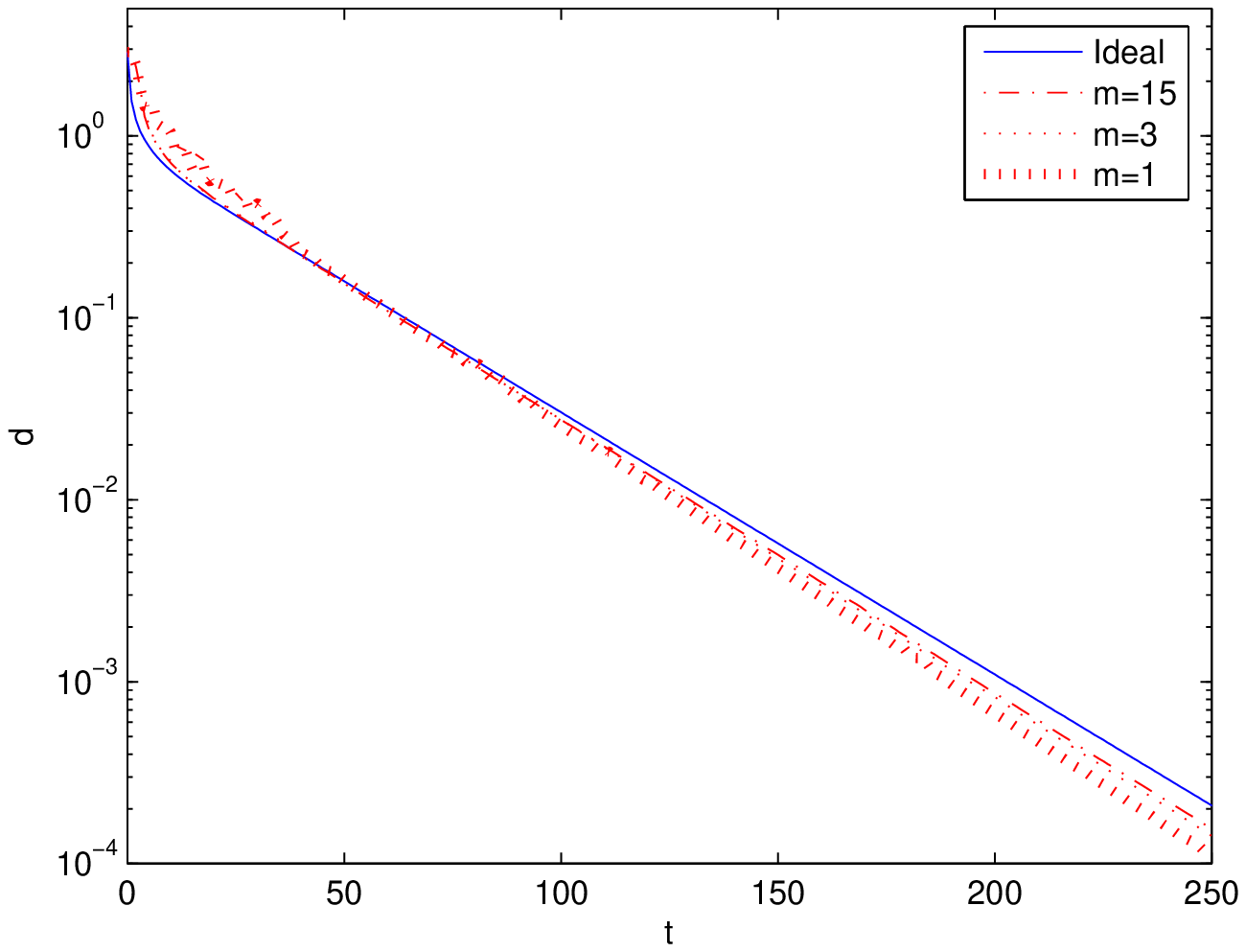}
\caption{Ring graph, in the $\kin<\rho$ regime. $N=20$, $\kin=0.9$, $\kout=2$.}\label{fig:ZoomRing_m2}
\end{figure}


\begin{figure}[htb]
\includegraphics[width=.8\textwidth]{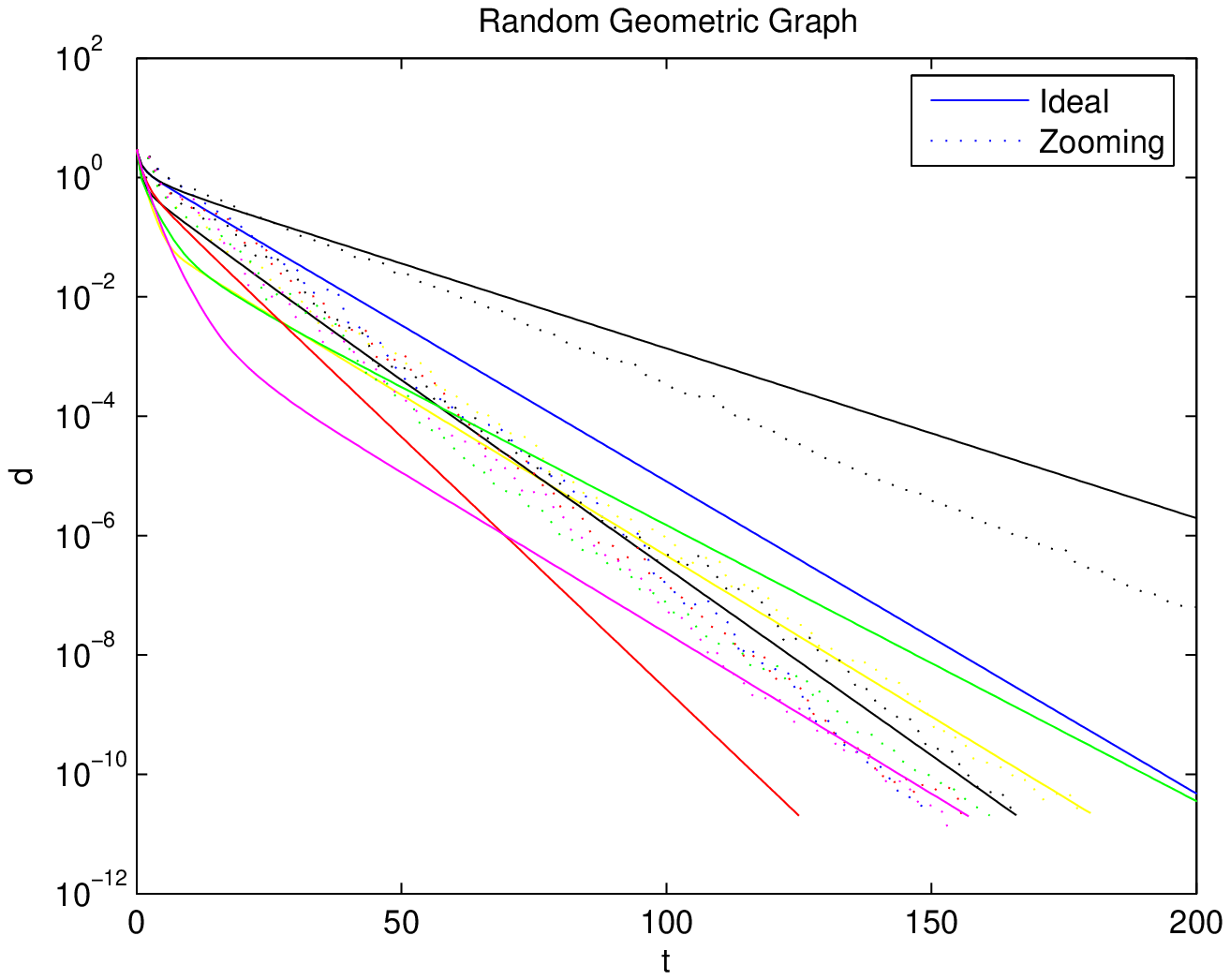}
\caption{Different performance on different samples of random geometric graph with $R=0.5$.  $m=3$, $\kin=0.3$, $\kout=2$.}\label{fig:ZoomRGGinst}
\end{figure}

Not to restrict ourselves to ring graphs, we show some simulations regarding random geometric graphs in Figure~\ref{fig:ZoomRGGinst}. The plots show how the performance depends on the topology of the geometric graph: different samples have different spectral radiuses, and then different convergence rates in the ideal case. Instead, in most cases the zooming algorithm is in the $\kin<\rho$ regime, and then the performance is roughly independent of the sample. This suggests that in this regime the algorithm is little sensitive to the spectral radius, and rises the question of which is a significant graph theoretic index.


\section{Conclusions}\label{sec:outro}
With theoretical and experimental results we showed that using the zooming scheme the average consensus problem can be efficiently solved although the agents can send only quantized information. Indeed the systems converge to average consensus at exponential speed, with a convergence rate which can be chosen as small as the convergence rate of consensus with perfect exchange of information, provided the number of quantization levels is large enough.
Though the theoretical results are quite conservative the efficiency of these methods is apparent from simulations, which show the method to converge outside the scope of current results, and converge at a rate faster than the ideal algorithm. Improving the theoretical analysis is thus a major concern, especially about the issue of scalability in the number of agents.

We have already underlined in Remark~\ref{rem:NeedForSynch} that the algorithm requires the states estimates to be shared among neighbors. This implies that the algorithm is not suitable, in the present form, for implementation on a digital noisy channel, unless we allow a feedback on the channel state. A natural development of this research is thus designing algorithms to solve the average consensus problem at exponential speed over a network of digital noisy channels.

\bibliographystyle{abbrv} \bibliography{aliasFrasca,PF,RefFrasca}
\end{document}